\documentclass[11pt,twoside]{article}

\usepackage{amsfonts}
\usepackage{amsmath}
\usepackage{amssymb}
\usepackage{amsthm}
\usepackage{graphics} 
\usepackage{epsfig} 
\usepackage{latexsym}

\graphicspath{{./figures/}} 

\usepackage[textwidth=16cm,textheight=22cm,top=3cm,bottom=3cm]{geometry}
\usepackage{fancyhdr}

\newtheorem{theorem}{Theorem}
\newtheorem{lemma}{Lemma}
\newtheorem{corollary}{Corollary}

\newtheorem{remark}{Remark}
\theoremstyle{definition}

\newtheorem{manualtheoreminner}{Theorem}
\newenvironment{manualtheorem}[1]{%
  \renewcommand\themanualtheoreminner{#1}%
  \manualtheoreminner
}{\endmanualtheoreminner}

\usepackage{color}

\newcommand{\drop}[1]{}
\newcommand{\no}{\noindent}
\newcommand{\fer}[1]{(\ref{#1})}
\newcommand{\qtext}[1]{\quad\text{#1}}
\newcommand{\qtextq}[1]{\quad\text{#1}\quad}

\newcommand{\ba}{\mathbf{a}}

\newcommand{\bn}{\mathbf{n}}

\newcommand{\bu}{\mathbf{u}}
\newcommand{\bv}{\mathbf{v}}

\newcommand{\cT}{\mathcal{T}}

\newcommand{\cB}{\mathcal {B}}

\newcommand{\eps}{\varepsilon}
\newcommand{\Deltan}{\Delta^{\!n}}
\newcommand{\vfi}{\varphi}

\newcommand{\grad}{\nabla}
\renewcommand{\o}{\omega}
\newcommand{\p}{\partial}
\newcommand{\chil}{\raisebox{0.5mm}{\large $\chi$}}

\newcommand{\N}{\mathbb{N}}

\newcommand{\R}{\mathbb{R}}

\def\O{\Omega}

\newcommand{\abs}[1]{\lvert #1 \rvert}
\newcommand{\nor}[1]{\| #1 \|}

\DeclareMathOperator{\Div}{div}

\DeclareMathOperator{\supp}{supp}
\DeclareMathOperator{\tr}{tr}
\DeclareMathOperator{\dist}{dist}
\DeclareMathOperator{\interior}{int}

\newcommand{\wto}{\rightharpoonup}

\title{Convergence of solutions of a rescaled evolution nonlocal cross-diffusion problem to its local diffusion counterpart
\thanks{Supported by Spanish MCI Project PID2020-116287GB-I00. 
}}

\author{Gonzalo Galiano, ~~~Juli\'an Velasco  \thanks{Dpt. of Mathematics, University of  Oviedo, Spain ({\tt galiano@uniovi.es, julian@uniovi.es})}}

\date{}

\begin{document}

\maketitle

\begin{abstract}
We prove that, under a suitable rescaling of the integrable kernel defining the nonlocal diffusion terms, the corresponding sequence of solutions of the Shigesada-Kawasaki-Teramoto  nonlocal cross-diffusion problem converges to a solution of the usual problem with local diffusion. In particular, the result may be regarded as a new proof of existence of solutions for the local diffusion problem.  
\vspace{0.25cm}

\no\emph{Keywords: }Nonlocal diffusion, integrable kernel, cross-diffusion, rescaled problem, convergence, Shigesada-Kawasaki-Teramoto population model.

\end{abstract}

\section{Introduction}

Evolution nonlocal diffusion escalar problems with integrable kernels have been extensively investigated in recent years, see the monograph by Andreu et al. \cite{Andreu2010} and the references therein. The paradigmatic problem is the evolution nonlocal $p-$Laplacian, which is expressed as: given $T>0$ and $\O\subset\R^N$ $(N\geq 1)$ an open set, find $u:[0,T]\times\O\to\R_+$ such that 
\begin{align}
&  \p_t u(t,x)  =  \int_\O J(x-y) \abs{u(t,y))-u(t,x)}^{p-2}\big(u(t,y)-u(t,x) \big) dy, \label{eq.eqp}\\
&   u(0,x)=  u_{0}(x),  \label{eq.idp}
\end{align}
 for $(t,x)\in Q_T=(0,T)\times\O$, and for some $u_{0}:\O\to\R_+$. Here, $\R_+=[0,\infty)$,  and the diffusion kernel, $J:\R^N\to\R_+$, is usually assumed to be continuous, radial, radially non-increasing, and with unitary norm in $L^1(\R^N)$.
 
 Terming the equation \fer{eq.eqp} as the \emph{evolution nonlocal $p-$Laplacian} equation has its reasons. Firstly, because it arises as the gradient descent of the Euler-Lagrange equation of the energy functional $E^{\text{nl}}_p(v)=\int_\O\int_\O J(x-y)\abs{v(y))-v(x)}^{p}dxdy$, in analogy to the usual evolution $p-$Laplacian equation, for which the energy is given by $E_p(v)=\int_\O\abs{\grad v}^{p}$. Secondly, and most importantly, because under the  rescaling $J_n(x)=n^{N+p}J(nx)$, the corresponding sequence of solutions, $u_n$, of problem \fer{eq.eqp}-\fer{eq.idp} converges to the usual weak solution of the local  $p-$Laplacian problem 
 \begin{align*}
  \p_t v  =   \Div\big(\abs{\grad v}^{p-2}\grad v\big), 
\end{align*}
with the same initial datum $v(0,\cdot)=  u_{0}$ and with homogeneous Neumann boundary conditions. 

The main idea behind the choice of the rescaling is that, being $n^{N}J(nx)$ an approximation of the Dirac delta, the factor $n^p$ plays the role of the denominator of the continuous incremental spatial ratio
\begin{align*}
	\Big\lvert\frac{v(t,x+\frac{1}{n}z)-v(t,x)}{1/n}\Big\rvert^p,
\end{align*}
so that one expects that, for all $w\in W^{1,p}(\O)$, 
\begin{align*}
	E^{\text{nl}}_{p,n}(w) & =\int\int J_n(x-y)\abs{w(y))-w(x)}^{p}dxdy \\
	&= \int\int J(z)\Big\lvert\frac{w(x+\frac{1}{n}z)-w(x)}{1/n}\Big\rvert^p dzdx 
	 \to \int\abs{\grad w(x)}^{p}dx = E_p(w) , 
\end{align*}
as $n\to\infty$. 

In fact, Andreu et al. \cite[Theorem~6.12]{Andreu2010} show, among other properties, that the sequence $\{u_n\}$ of solutions of the nonlocal rescaled problems converge to $v$, the solution of the local problem, strongly in $L^\infty(0,T;L^p(\O)$. There are two main ingredients in their proof. The first is the precompactness result \cite[Theorem~6.11]{Andreu2010} based on previous results by Bourgain et al. \cite[Theorem~4]{Bourgain2001}, which shows that if 
$w_n \wto w$ weakly in $L^p(Q_T)$ and $E^{\text{nl}}_{p,n}(w_n)$ is uniformly bounded then $w_n\to w$ strongly in $L^p(\O)$ and $w\in W^{1,p}(\O)$. The second ingredient is the monotonocity of $\abs{s}^{p-2}s$, which plays an important role both in the theory of existence of solutions and in the identification of the limit of the solutions of the rescaled problems.

The objective of this article is to show that the same convergence property is true for a class of evolution nonlocal cross-diffusion systems. In \cite{Galiano2019b} (see also \cite{Galiano2019a} for related work), we introduced and proved the existence of solutions of the nonlocal version of the paradigmatic  Shigesada-Kawasaki-Teramoto (SKT) population model \cite{Shigesada1979}. The nonlocal version of this model is the following: for $i=1,2$, find $u_i:[0,T]\times\O\to\R_+$ such that 
\begin{align}
&  \p_t u_i(t,x)  =  \int_\O J(x-y) \big(p_i (\bu(t,y))-p_i (\bu(t,x)) \big) dy +f_i(\bu(t,x)), \label{eq.eq}\\
&   u_i(0,x)=  u_{0i}(x),  \label{eq.id}
\end{align}
 where  $\bu=(u_1,u_2)$ and, for $i,j=1,2$, $i\neq j$, the diffusion and reaction functions are given by
  \begin{align*}
 p_i(\bu) = u_i(c_i+a_i u_i + u_j),\quad  f_i(\bu) = u_i\big(\alpha_i-(\beta_{i1} u_1 + \beta_{i2}u_2)\big),
\end{align*}
for some non-negative constant coefficients $c_i,~a_i,~\alpha_i,~\beta_{ij}$.

The local diffusion problem, i.e., the original version of the SKT model, is to find, for $i=1,2$, functions $v_i:[0,T]\times\O\to\R_+$ such that
\begin{align}
&  \p_t v_i  =  \Delta p_i (\bv) +f_i(\bv)&&\text{in }Q_T, & \label{eq.eqs}\\
&  \grad p_i (\bv) \cdot \bn = 0 &&\text{on }(0,T)\times\p\O, & \label{eq.bcs}\\
&   v_i(0,\cdot)=  u_{0i}  && \text{in }\O. \label{eq.ids}
\end{align}
The existence of solutions of evolution cross-diffusion problems like \fer{eq.eqs}-\fer{eq.ids} has been addressed for a variety of problems \cite{Amann1989, Galiano2003, Chen2004, Galiano2014, Desvillettes2014, Jungel2015}, for which one can define an appropriate Lyapunov functional, also known as \emph{entropy functional}, which in the case of the SKT model is given by
\begin{align}
\label{def.ent}
 E_{\bv}(t) = \sum_{i=1}^2 \int_\O \big(v_i(t,\cdot)(\ln(v_i(t,\cdot))-1)+1) \geq0.
\end{align}
Formally testing  \fer{eq.eqs} with  $\ln(v_i)$ yields the estimate  
\begin{align}
\label{est.entropy}
 E_{\bv}(t)+ \sum_{i=1}^2 a_i \int_{Q_t} \abs{\grad v_i}^2 \leq E_{\bv}(0)+C,
\end{align}
which is the first step for a compactness argument for a sequence of solutions of regularized problems. The second step, with the horizon of applying Aubin-Lion's type compactness arguments, is to deduce suitable uniform estimates for the time derivatives $\p_t v_i$. For the local SKT problem, the fundamental tool to obtain these estimates is the Gagliardo-Niremberg inequality, although other approaches (for other cross-diffusion problems), like global regularity in H\"older spaces \cite{Amann1989} or duality estimates \cite{Desvillettes2014}, are also fruitful.

For the nonlocal SKT problem a nonlocal entropy estimate analogous to \fer{est.entropy} is deduced by similar procedures than in the local case, obtaining an estimate which is independent of the kernel $J$, see \fer{th1:entropy}.  A fundamental difference between the local and the nonlocal models is that, in the latter, the $L^\infty(Q_T)$ regularity of solutions is proven, although with a bound depending on $J$. As a consequence, the time derivative bound is simpler to achieve than in the local case, and moreover, other important properties, like the uniqueness of solutions, are deduced based on this regularity.

Turning to the rescaled problems, the entropy estimate is again deduced as in the proof of existence of solutions, since this estimate is independent of the kernel. However, the estimation of the time derivative is more problematic. In the proof of our main theorem, and in view of the lack of a suitable nonlocal Gagliardo-Nirenberg inequality, we have resorted to the use of duality estimates, see e.g. \cite{Pierre1997, Desvillettes2014}, which allow us to show an improved regularity of the sequence of rescaled solutions (with respect to that implied by the entropy estimate).  

Finally, to close this introduction, let us notice that an interesting by-product of the main result of this article is to regard it as a new proof of the existence of solutions of the usual local diffusion SKT problem. Besides, since the solutions of the approximating problems do not need to have  Sobolev regularity, our approach opens the possibility of approximating numerically the solutions of the local SKT problem by alternative methods \cite{Perez2011, Du2019, Delia2020A, Delia2020B,  Galiano2021}.

\section{Hypothesis and main result }

We always assume, at least, the following hypothesis on the data, that we shall refer to as \textbf{(H)}:
\begin{enumerate}
\item The final time, $T>0$, is arbitrarily fixed. The spatial domain, $\O\subset\R^N$ $(N\geq 1)$, is an open, bounded and Lipschitz set.

\item The kernel function $J:\R^N\to\R$ is a non-negative, continuous, radial, radially non-increasing function with compact support and such that 
  \begin{align*}
  \int_{\R^N} J(x)d x = 1.
  \end{align*}

\item The initial data $u_{0i}\in L^\infty(\O)\cap BV(\O)$ are non-negative, for $i=1,2$. 

\item For $i,j=1,2$, the constants $c_i,~a_i,~\alpha_i,~\beta_{ij}$ are non-negative. 
\end{enumerate}  
We recall here the result on existence and uniqueness of solution of problem \fer{eq.eq}-\fer{eq.id} proven in  \cite{Galiano2019b}. Notice that Hypothesis (H)$_2$, needed for the result on convergence of the rescaled problems, is more restrictive than the corresponding assumption in \cite{Galiano2019b} used to prove the existence and uniqueness of solution. We merged both conditions for the sake of brevity. In particular, this assumption implies the existence of some $r>0$ such that 
\begin{align}
\label{suppJ}
\supp(J)=B_r(0).	
\end{align}
Regarding the restriction $u_{0i}\in BV(\O)$, unusual for the local diffusion problem, we must include it since it is needed for the compactness argument employed in the proof of the existence of solutions of the nonlocal diffusion model.

\begin{manualtheorem}{A}[\textbf{Existence and uniqueness of solution \cite{Galiano2019b}}] \label{th.existence}
Assume (H) and
\begin{align*}
 a_i+\beta_{ii} >0, \qtext{for } i=1,2 .
\end{align*}
Then, there exists a unique strong solution $(u_1,u_2)$ of problem \fer{eq.eq}-\fer{eq.id} with $u_i\geq 0$ a.e. in $Q_T$ and such that, 
for $i=1,2$ and $t\in [0,T]$,
\begin{align}
&  u_i\in W^{1,\infty}(0,T;L^\infty(\O))\cap C([0,T];L^\infty(\O)\cap BV(\O)), \nonumber\\
& E_{\bu}(t) + \sum_{i=1}^2 a_i \int_{Q_t}\int_\O J(x-y) \big(u_i(s,x)-u_i(s,y)\big)^2 dy dx ds \leq E_{\bu}(0)+C, \label{th1:entropy}
\end{align}
with $E_{\bu}(t)$ defined by \fer{def.ent}, and for some constant $C>0$ independent of $J$.
\end{manualtheorem}
Although not mentioned in the above theorem, integrating \fer{eq.eq} in $(0,t)\times\O$ yields the  uniform estimate 
\begin{align}
\label{est.L1}
\nor{u_{i}}_{L^\infty(0,T;L^1(\O))} \leq C\nor{u_{0i}}_{L^1(\O)},
\end{align}
where, here and in what follows, $C$ denotes a generic positive constant independent of $J$.

For stating the main result contained in this article, we introduce the rescaled kernel
 \begin{align}
 \label{def.jeps}
  J_n(x) = C_1 n^{2+N}J(nx),\qtext{with } C_1^{-1} = \frac{1}{2} \int_{\R^N} J(x)x_N^2 dx.
 \end{align}
The kernel $J_n$ satisfies the conditions of Theorem~\ref{th.existence} and, therefore,  
there exists a unique solution, $\bu_n$, of \fer{eq.eq}-\fer{eq.id}, corresponding to $J_n$. We also introduce here, for later reference, the approximation to the Dirac delta
\begin{align}
 \label{def.jtilde}
  \tilde J_n(x) = n^N J(nx) = \frac{1}{C_1 n^2}J_n(x).	
\end{align}

\begin{theorem}
\label{th:convergence}
	Assume (H) and suppose that $a_i>0$ for $i=1,2$. Let  $q=3-\delta$, with $\delta>0$ a small number, and $r=2q/(q+2)$. Let $\bu_n$ be the solution of problem \fer{eq.eq}-\fer{eq.id} corresponding to the kernel $J_n$ and consider the sequence $\{\bu_n\}$. There exists a subsequence $\{\bu_{n_j}\}$ such that $\bu_{n_j}\to \bu$ strongly in $L^q(Q_T)$, where $\bu\in W^{1,6/5}(0,T;(W^{1,6}(\O))')\cap L^2(0,T;H^1(\O))\cap L^q(Q_T)$ is a weak solution of problem \fer{eq.eqs}-\fer{eq.ids} in the following weak sense:
\begin{align*}
  \int_0^T <\p_t u_{i}, \xi>    + \int_{Q_T} \grad p_i(\bu) \cdot \grad \xi  = \int_{Q_T} f_i(\bu) \xi,
\end{align*}
for all  $\xi\in L^{r'}(0,T;W^{1,r'}(\O))$, 
where $<\cdot,\cdot>$ denotes the duality product of $(W^{1,r'}(\O))'\times W^{1,r'}(\O)$, and 
\begin{align}
\label{idconv}
  \int_0^T <\p_t u_{i}, \psi>    + \int_{Q_T} (u_i-u_{i0})\p_t \psi = 0,
\end{align}
for all $\psi \in L^{6}(0,T;W^{1,6}(\O)) \cap W^{1,q'}(0,T;L^{q'}(\O))$. In addition, $\bu$ satisfies the entropy estimate \fer{est.entropy}.
\end{theorem}
Notice that the regularity implied by the definition of the exponents $q$ and $p$ is independent of the spatial dimension, $N$. Regarding Theorem~\ref{th:convergence} as a proof of existence of solutions of problem \fer{eq.eqp}-\fer{eq.idp},  our result improves that of \cite{Chen2004} for $N\geq 3$. For $N=1,2$, the regularity proven in \cite{Galiano2003, Chen2004} may be 
also recovered in our solution by using the Sobolev imbbeding theorem ($N=1$) or the Gagliardo-Niremberg inequality ($N=2$).

Due to the lack of knowledge concerning the uniqueness of solution of the local diffusion problem  \fer{eq.eqp}-\fer{eq.idp}, the convergence of the full sequence of rescaled problems may not be ensured, but only that of a subsequence to some solution of \fer{eq.eqp}-\fer{eq.idp}.

\section{Proof of Theorem~\ref{th:convergence} }
We start recalling a fundamental result by Bourgain et al. \cite[Theorem~4]{Bourgain2001}. We use the variant introduced by Andreu et al. \cite[Theorem~6.11]{Andreu2010} and state a straightforward extension to deal with time-dependent functions. First, we introduce the notation for the extension by zero of a function $\psi:\O\to\R$:
\begin{align*}
\bar\psi(x) = \begin{cases}
 \psi(x) & \text{if } x\in \O,	\\
 0 & \text{if } x\in \O^c.
 \end{cases}
\end{align*}

\begin{manualtheorem}{B}[\cite{Bourgain2001, Andreu2010}]
\label{th.brezis}	
Let $1\leq p <\infty$ and assume that $\O\subset\R^N$ satisfies (H)$_1$. Let $\rho:\R^N\to\R$ be a non-negative, continuous, radial, radially non-increasing function with compact support, and set $\rho_n(x)=n^N\rho(nx)$.
Let $\{f_n\}$ and $\{g_n\}$ be sequences in $L^p(\O)$ and $L^p(Q_T)$, respectively, such that 
\begin{align}
& \int_{\O}\int_\O \rho_n(x-y) \big(f_{n}(x)-f_{n}(y)\big)^p dy dx \leq \frac{C_0}{n^p}, \label{BQ1a}\\
&   \int_{Q_T}\int_\O \rho_n(x-y) \big(g_{n}(t,x)-g_{n}(t,y)\big)^p dy dx dt \leq \frac{C_0}{n^p}. \nonumber 
\end{align}
Let $\Phi_\delta = \frac{1}{\abs{B_\delta(0)}}\chil_{B_\delta(0)}$, for $\delta>0$. Then, there exists a constant  $C\equiv C(p,\O,\rho)$ independent of $n$,  and a number $n_\delta\in\N$ such that, for $n\geq n_\delta$, 

\no (a)
\begin{align}
&\int_{\O} \abs{f_n}^p  \leq C\Big(C_0 + \Big\lvert\int_\O f_n\Big\rvert ^p\Big),\label{brezis.lq}\\
&\int_{\O}  \abs{f_n-f_n*\Phi_\delta}^p \leq CC_0 \delta^p 
.	\label{brezis.wq}
\end{align}
In consequence, there exists a subsequence $\{f_{n_k}\}$ and a function $f\in W^{1,p}(\O)$ ($BV(\O)$, if $p=1$) such that 
\begin{align}
\label{brezis.convq}
f_{n_k} \to f \qtextq{strongly in} L^p(\O).	
\end{align}

\no (b)
\begin{align}
&\int_{Q_T} \abs{g_n}^p  \leq C\Big(C_0 + \Big\lvert\int_{Q_T} g_n\Big\rvert ^p\Big),\label{brezis.lqT}\\
&\int_{Q_T}  \abs{g_n(t,x)-(g_n(t,\cdot)*\Phi_\delta)(x)}^p dxdt\leq CC_0\delta^p . \nonumber	
\end{align}
In addition, if $g_n\wto g$ weakly in $L^p(Q_T)$ for $1<p<\infty$, then 
\begin{align}
\label{lemma.conv2}
\rho(z)^{1/p}\chil_\O \Big(x+\frac{z}{n}\Big)	\frac{\bar g_n(t,x+\frac{1}{n}z)-g_n(t,x)}{1/n} \wto \rho(z)^{1/p} ~ z\cdot \grad g(t,x)
\end{align}
weakly in $L^p(Q_T)\times L^p(0,T;L^p(\R^N))$.
\end{manualtheorem}
\no\emph{Proof. }Estimates \fer{brezis.lq} and \fer{brezis.wq}, as well as their consequence \fer{brezis.convq}, are proven in \cite{Bourgain2001, Andreu2010}. Their extension to the time-dependent functions of (b) is straightforward once we notice that the results of (a) are valid pointwise for a.e. $t\in (0,T)$. Finally, \fer{lemma.conv2} is also a direct consequence of \cite[Theorem~6.11 (1)]{Andreu2010}. $\Box$

\begin{remark}
The inequality \fer{brezis.lq} is a Poincar\'e's type inequality, since the constant $C_0$ may be replaced by the nonlocal energy appearing in \fer{BQ1a}, see the proof of \cite[Theorem~4]{Bourgain2001} 	for details. Thus, this provides an alternative (constructive) proof of the result stated in \cite[Proposition~6.19]{Andreu2010}. 
\end{remark}

\no\textbf{Step 1. Uniform bound in $L^3(Q_T)$. }The entropy inequality \fer{th1:entropy} and the $L^1(Q_T)$ estimate \fer{est.L1} together with \fer{brezis.lqT} applied with $p=2$, $\rho_n = \tilde J_n$, and $g_n=u_{in}$ imply that the sequences $\{u_{in}\}$, for $i=1,2$, are uniformly bounded in $L^2(Q_T)$. However, this bound is not enough to define the weak limit in an appropriate  reflexive space of test functions since the nonlinear parts of the limit diffusive term are expected to be of the form
\begin{align*}
\int_{Q_T} u_j \grad u_i\cdot\grad\vfi,	
\end{align*}
with the regularity  $\grad u_i \in L^2(Q_T)$. Thus, we start improving the uniform bounds of $\{u_{in}\}$ to the space  $L^3(Q_T)$.
This bound is obtained using an estimate of the dual problem corresponding to the nonlocal rescaled problem. The existence of solutions of this dual problem is ensured by the following lemma. 
\begin{lemma}
\label{lemma:dual}
	The problem: find $\phi:Q_T\to\R$ such that 
\begin{align}
&\p_t \phi (t,x) - a(t,x)\int_\O \rho(x-y) (\phi(t,y)-\phi(t,x))dy = b(t,x), \label{eq.test1}\\
&\phi(0,x) = c(x),	\label{eq.test2}
\end{align}
where $\rho \in L^\infty(\R^N)$, $a,b \in L^\infty(Q_T)$, and $c\in L^\infty(\O)$, has a unique solution such that 
\begin{align}
\label{test.reg}
\phi, ~\Delta^{\!1,\rho} \phi \in C([0,T];L^\infty(\O)),	 \text{ and } \p_t \phi \in L^\infty(Q_T), 
\end{align}
where $\Delta^{\!1,\rho} \phi(t,x)=\int_\O \rho(x-y) (\phi(t,y)-\phi(t,x))dy$.
In addition, if $a,b$ and $c$ are non-negative then $\phi$ is non-negative. 
\end{lemma}
The proof of this and the other lemmas used for proving Theorem~\ref{th:convergence} are given in Section~\ref{sec:lemmas}. 
In the following, we shall use the notation
\begin{align*}
&	\Deltan \vfi(t,x) = \int_\O J_n(x-y)(\vfi(t,y)-\vfi(t,x)) dy, \qtextq{for any} \vfi\in L^1(Q_T),\\
&	\tilde p_i(\bu) := p_i(\bu)/u_{i}= \alpha_i -(\beta_{i1}u_1+\beta_{i2u_2}).
\end{align*}

\begin{corollary}
\label{cor.test}
	Let $\vfi_{in} (t,x)= e^{-\lambda(T-t)}\phi(T-t,x)$, where $\lambda>0$ is a constant and $\phi$ is the non-negative solution of \fer{eq.test1}-\fer{eq.test2} corresponding to $\rho(x)=J_n(x)$, $a(t,x)=\tilde p_i(\bu_n(T-t,x))$, $b(t,x)=-e^{\lambda t}\psi(T-t,x)\sqrt{\tilde p_i(\bu_n(T-t,x))}$, and $c=0$, being $\psi\in L^\infty(Q_T)$  a non-positive arbitrary function. Then $\vfi_{in}$  is a non-negative solution of 
\begin{align}
&\p_t\vfi_{in}+ \tilde p_i(\bu_n) \Deltan \vfi_{in} - \lambda  \vfi_{in} = \sqrt{\tilde p_i(\bu_n)} \psi& & \text{in }Q_T, &\label{eq.test3}\\
&\vfi_{in}(T,\cdot) = 0 & & \text{on }\p\O, &	\label{eq.test4}
\end{align}
with the same regularity than $\phi$, see \fer{test.reg}.
\end{corollary}
\no\emph{Proof of Corollary~\ref{cor.test}.}
By Theorem~\ref{th.existence}, we have that $\tilde p_i(\bu_n) \in L^\infty(Q_T)$  is non-negative, so that $\rho,a,b,c \in L^\infty$ are non-negative. By Lemma~\ref{lemma:dual}, there exists a unique non-negative solution  $\phi$ of problem \fer{eq.test1}-\fer{eq.test2} corresponding to this data. A simple calculation shows that $\vfi_{in}$ is then the non-negative solution of  \fer{eq.test3}-\fer{eq.test4} inheriting the same regularity than $\phi$. $\Box$

\medskip

Now we proceed to obtain the $L^3(Q_T)$ uniform bound of the sequences $\{u_{in}\}$, for $i=1,2$. We multiply the equation \fer{eq.eq} of $u_{i n}$ by the solution $\vfi_{in}$ of problem \fer{eq.test3}-\fer{eq.test4} and integrate to get, for $i=1,2$,
\begin{align*}
	\int_\O u_{i0}\vfi_{in}(0,\cdot) &+\int_{Q_T} u_{i n}\Big(\p_t\vfi_{in}+ \tilde p_i(\bu_n) \Deltan \vfi_{in} - \lambda  \vfi_{in}\Big) \\
	&=  -\int_{Q_T}\big(\lambda u_{in}+f_i(\bu_n)\big)\vfi_{in} . 
\end{align*}
Using the equation \fer{eq.test3},  the explicit expression of $f_i$ and the non-negativity of $u_{in}$ and $\vfi_{in}$, we obtain 
\begin{align*}
-\int_{Q_T} u_{i n}\sqrt{\tilde p_i(\bu_n)} \psi \leq \int_\O u_{i0}\vfi_{in}(0,\cdot) +(\alpha_i+\lambda)\int_{Q_T}u_{in}\vfi_{in}.
\end{align*}
Noticing that $\psi\leq 0$ and using H\"older's inequality and the uniform estimate of $u_{in}$ in $L^2(Q_T)$, see Step 1,  we obtain 
\begin{align}
\label{dual3}
\Big\lvert\int_{Q_T} u_{i n}\sqrt{\tilde p_i(\bu_n)} \psi \Big\rvert \leq   C\big(\nor{\vfi_{in}(0,\cdot)}_{L^{2}(\O)} +\nor{\vfi_{in}}_{L^{2}(Q_T)} \big).
\end{align} 
Our objective is to estimate $\nor{\vfi_{in}(0,\cdot)}_{L^{2}(\O)}$ and $\nor{\vfi_{in}}_{L^{2}(Q_T)}$ in terms of $\nor{\psi}_{L^2(Q_T)}$ to deduce, by duality, a uniform estimate of $\nor{u_{i n}\sqrt{\tilde p_i(\bu_n)}}_{L^2(Q_T)}$.

Multiplying the equation \fer{eq.test3} of $\vfi_{in}$ by $\Deltan \vfi_{in}$ and using the nonlocal integration by parts formula and Young's inequality yields
\begin{align}
	\int_{Q_T} \p_t\vfi_{in} \Deltan \vfi_{in} &+ \int_{Q_T} \tilde p_i(\bu_n) \abs{\Deltan \vfi_{in}}^2\nonumber \\
	&+ \frac{\lambda}{2}\int_{Q_T}\int_\O J_n(x-y)(\vfi_{in}(t,y)-\vfi_{in}(t,x))^2dydxdt \nonumber \\
	& \leq  \frac{1}{2}\int_{Q_T} \tilde p_i(\bu_n) \abs{\Deltan \vfi_{in}}^2 + \frac{1}{2}\int_{Q_T} \psi^2 .\label{dual1}
\end{align} 
We have
\begin{align*}
	\int_{Q_T} \p_t\vfi_{in} \Deltan \vfi_{in} &= \int_{Q_T} \p_t \vfi_{in}(t,x)  \int_\O J_n(x-y)(\vfi_{in}(t,y)-  \vfi_{in}(t,x)) dydxdt  \\
    &=	-\frac{1}{4}\int_{Q_T} \p_t \int_\O  J_n(x-y)(\vfi_{in}(t,y)-  \vfi_{in}(t,x))^2  dydxdt\\
    &=	\frac{1}{4}\int_{\O} \int_\O  J_n(x-y)(\vfi_{in}(0,y)-  \vfi_{in}(0,x))^2  dydx    ,
\end{align*}
so that from \fer{dual1}, we obtain
\begin{align}
	\frac{1}{4}\int_{\O} \int_\O  J_n(x-y) & (\vfi_{in}(0,y)-  \vfi_{in}(0,x))^2  dydx + \frac{1}{2}\int_{Q_T} \bar p_i(\bu_n) \abs{\Deltan \vfi_{in}}^2\nonumber \\
	& + \frac{\lambda}{2}\int_{Q_T}\int_\O J_n(x-y)(\vfi_{in}(t,y)-\vfi_{in}(t,x))^2dydxdt\nonumber \\
&	\leq   \frac{1}{2}\int_{Q_T} \psi^2. \label{dual2}
\end{align} 
Therefore, using Theorem~\ref{th.brezis} (a) with $f_n=\vfi_{in}(0,\cdot)$ and (b) with $g_n=\vfi_{in}$ we deduce from \fer{dual2}
\begin{align}
& 	\nor{\vfi_{in}(0,\cdot)}_{L^2(\O)}^2\leq C\Big(\nor{\psi}_{L^2(Q_T)}^2 + \Big\lvert\int_\O \vfi_{in}(0,\cdot)\Big\rvert ^2\Big) , 	\label{dual9}\\
&	\nor{\vfi_{in}}_{L^2(Q_T)}^2\leq C\Big(\nor{\psi}_{L^2(Q_T)}^2 + \Big\lvert\int_{Q_T} \vfi_{in}\Big\rvert ^2\Big)  . 	\label{dual8}
\end{align}
 Integrating the equation \fer{eq.test3} of $\vfi_{in}$ in $(t,T)\times \O$ yields   
\begin{align*}
	\int_\O \vfi_{in}(t,\cdot) 
	& = \int_t^T \int_\O \tilde p_i(\bu_n))  \Deltan \vfi_{in} 
	- \lambda \int_t^T \int_\O \vfi_{in}
	 - \int_t^T \int_\O \sqrt{\tilde p_i(\bu_n)} \psi	 .
\end{align*}
Observing that $\vfi_{in}$ and $\lambda$ are non-negative and using H\"older's inequality, we get
\begin{align*}
	\int_\O \vfi_{in}(t,\cdot) 
	& \leq \Big(\int_{Q_T} \tilde p_i(\bu_n) \Big)^{1/2}\Big[ \Big(\int_{Q_T} \tilde p_i(\bu_n) \abs{\Deltan \vfi_{in}}^2 \Big)^{1/2}  + \Big(\int_{Q_T} \psi^2 \Big)^{1/2}\Big],
\end{align*}
and, therefore, 
\begin{align*}
&	\int_\O \vfi_{in}(0,\cdot) \leq 2 \Big(\int_{Q_T} \psi^2 \Big)^{1/2}\Big(\int_{Q_T} \tilde p_i(\bu_n) \Big)^{1/2},\\
&	\int_{Q_T} \vfi_{in} \leq 2T \Big(\int_{Q_T} \psi^2 \Big)^{1/2}\Big(\int_{Q_T} \tilde p_i(\bu_n) \Big)^{1/2}.
\end{align*}
Returning to \fer{dual9}-\fer{dual8} and taking into account the uniform $L^1(Q_T)$ estimate \fer{est.L1}, we deduce 
\begin{align*}
 	\nor{\vfi_{in}(0,\cdot)}_{L^2(\O)}\leq C\nor{\psi}_{L^2(Q_T)},\quad  
	\nor{\vfi_{in}}_{L^2(Q_T)}\leq C\nor{\psi}_{L^2(Q_T)}. 
\end{align*} 
Finally,  \fer{dual3} yields, by duality, an uniform estimate for
	$\nor{u_{i n}\sqrt{\tilde p_i(\bu_n)}}_{L^{2}(Q_T)}$, or, in other words, the estimate
\begin{align*}
 \int_{Q_T} u_{i n}^2(c_i + a_iu_{i n} + u_{jn}) \leq C.
\end{align*}
In particular, if $a_i>0$, we obtain uniform estimates of $u_{i n}$ in $L^3(Q_T)$.

\medskip

\no\textbf{Step 2. Strong convergence. }
The following lemma is a consequence of two results concerning compactness: the precompactness result of Bourgain et al. \cite[Theorem~4]{Bourgain2001} for sequences defined on the spatial domain $\O$,  and  the compensated compactness result of P.L. Lions \cite[Lemma~5.1]{Lions} for the product of sequences defined in $Q_T$. 
\begin{lemma}
\label{lemma:compactness}
Let $\rho:\R^N\to\R$ be like in Theorem~\ref{th.brezis}, and let $\{f_n\}$ be a sequence in $L^3(Q_T)$ such that 
\begin{align}
&(i)~ \nor{f_n}_{L^3(Q_T)}\leq C,\label{BQ0}\\
   &(ii)~\int_{Q_T}\int_\O \rho_n(x-y) \big\lvert f_{n}(t,x)-f_{n}(t,y)\big\rvert^2 dy dx dt \leq \frac{C}{n^2}, \label{BQ1}\\
   &(iii)~\p_t f_n \text{ is uniformly bounded in }L^1(0,T;W^{-m,1}(\O)), \label{BQ2}
\end{align}
for some $m\geq 0$ independent of $n$.
Then there exists a subsequence $\{f_{n_k}\}$ and a function $f \in L^q(Q_T)\cap L^2(0,T;H^{1}(\O))$ such that
\begin{align}
\label{lemma.conv1}
f_{n_k}\to f	\text{ strongly in }L^q(Q_T) \qtextq{for any} q<3.
\end{align}
\end{lemma}

Taking $\rho = J$ and $f_n=u_{in}$ for $i=1,2$, the uniform estimate of $u_{in}$ in $L^3(Q_T)$ obtained in Step 1 and the entropy inequality satisfied by these functions, see \fer{th1:entropy}, imply (i) and (ii) of Lemma~\ref{lemma:compactness}. We now check that the uniform time estimate (iii) is also satisfied.  For any smooth function $\xi:Q_T\to\R$, we have
\begin{align}
  \int_{Q_T} \p_t u_{i n} \xi =&   \int_{Q_T} p_i (\bu_n(t,x)) \int_\O J_n(x-y) (\xi(t,y)-\xi(t,x))  dydxdt  \nonumber \\
  & + \int_{Q_T} f_i(\bu_{n}) \xi \leq c \big(\nor{p_i(\bu_n)}_{L^{3/2}(Q_T)}\nor{\Deltan \xi}_{L^3(Q_T)} \nonumber  \\
& +\nor{f_i(\bu_n)}_{L^{3/2}(Q_T)}\nor{\xi}_{L^3(Q_T)}\big). \label{est.time3}
\end{align}
 \begin{lemma}
 \label{w023}
 Let $\xi \in L^p(0,T;W^{2,p}_0(\O))$, for $1\leq p <\infty$. Then there exist a constant $C$ independent of $n$ and a constant $n_J\in\N$ such that if $n>n_J$ then
 \begin{align*}
\nor{\Deltan \xi}_{L^p(Q_T)} \leq C \nor{\xi}_{L^p(0,T;W^{2,p}_0(\O))}. 	
 \end{align*}
 \end{lemma}
Using this lemma with $p=3$ and noting that 
$u_{in}$ is uniformly bounded in $L^3(Q_T)$ we obtain  from \fer{est.time3}, by duality,
\begin{align}
\label{est.time}
\nor{\p_t u_{in}}_{L^{3/2}(0,T;W^{-2,3/2}(\O))} \leq C.	
\end{align}
Therefore, (iii) of Lemma~\ref{lemma:compactness} is satisfied and \fer{lemma.conv1} follows, this is, there exist  subsequences (not relabeled) such that $u_{in}\to u_i$ strongly in $L^q(Q_T)$, for $i=1,2$ and for any $1\leq q<3$. 

\medskip

\no\textbf{Step 3. Time derivative estimate. }
Once we have proven the strong convergence of $\{u_{in}\}$ in $L^q(Q_T)$, we may improve the uniform time estimate obtained in \fer{est.time}.  For $\xi$ smooth,  we have 
\begin{align}
  \int_{Q_T} \p_t u_{i n} \xi & =  -\frac{1}{2} \int_{Q_T}   \int_\O J_n(x-y)(p_i (\bu_n(t,y))-p_i (\bu_n(t,x)))   \nonumber \\
  & \times (\xi(t,y)-\xi(t,x))  dydxdt + \int_{Q_T} f_i(\bu_{n}) \xi = I_1+I_2 . \label{ts.id}
\end{align}
Clearly, $\abs{I_2}\leq C\nor{\xi}_{L^3(Q_T)}$, since $f_i$ is quadratic and $\{u_{in}\}$ is uniformly bounded in $L^{3}(Q_T)$. We  examine the terms of $p_i(\bu)=c_i u_i +a_iu_i^2+u_iu_j$ separately. For the linear term, we have, for $i=1,2$,
\begin{align*}
& \Big\lvert \int_{Q_T}   \int_\O J_n(x-y)(u_{in}(t,y)-u_{in}(t,x))  (\xi(t,y)-\xi(t,x))  dydxdt \Big\rvert\\
&  \leq \Big(\int_{Q_T} \int_\O J_n(x-y)\abs{u_{in}(t,y)-u_{in}(t,x)}^2 dydxdt\Big)^{1/2}\\
& \times
\Big(\int_{Q_T} \int_\O J_n(x-y)\abs{\xi(t,y))-\xi(t,x)}^2 dydxdt\Big)^{1/2} \leq C \nor{\xi}_{L^2(0,T;H^1(\O))},
\end{align*}
where we used the entropy estimate \fer{th1:entropy} and an straightforward modification of \cite[Theorem~1]{Bourgain2001} for including the time variable. For the quadratic terms, we have, for $i,j=1,2$, 
\begin{align}
& \Big\lvert \int_{Q_T}   \int_\O J_n(x-y)(u_{in}(t,y)u_{jn}(t,y)-u_{in}(t,x)u_{jn}(t,x)) \nonumber \\
& \times  (\xi(t,y)-\xi(t,x))  dydxdt \Big\rvert \leq I_{ji}+I_{ij} ,\label{decomp}
\end{align}
where, for $k,\ell = 1,2$, 
\begin{align*}
I_{k\ell} & = 	\Big\lvert \int_{Q_T}   \int_\O J_n(x-y)u_{kn}(t,y)(u_{\ell n}(t,y)-u_{\ell n}(t,x))  (\xi(t,y)-\xi(t,x))  dydxdt \Big\rvert .
\end{align*}
Since $J_n(z)= C_1n^2\tilde J_n(z)$, being $\tilde J_n$ an approximation of the Dirac delta, we have
\begin{align*}
I_{k\ell} \leq & C_1	 \Big(\int_{Q_T}   \int_\O \tilde J_n(x-y) \abs{u_{kn}(t,y)}^3 dydxdt\Big)^{1/3} \\
& \times \Big(\int_{Q_T}   \int_\O \tilde J_n(x-y)\frac{\abs{u_{\ell n}(t,y)-u_{\ell n}(t,x)}^2}{1/n^2} dydxdt\Big)^{1/2}\\
& \times \Big(\int_{Q_T}   \int_\O \tilde J_n(x-y)\frac{\abs{\xi(t,y)-\xi(t,x)}^{6} }{1/n^6}dydxdt\Big)^{1/6} .
\end{align*}
The first factor is bounded due to the uniform $L^3(Q_T)$ estimate found in Step~1. The second factor is bounded due to the entropy estimate  \fer{th1:entropy}. Finally, a new use of \cite[Theorem~1]{Bourgain2001} yields that the third factor is bounded by $\nor{\xi}_{L^6(0,T;W^{1,6}(\O))}$. This is, we obtain  
\begin{align*}
I_{k\ell} \leq C \nor{\xi}_{L^{6}(0,T;W^{1,6}(\O))}	.
\end{align*}
Thus, by duality, we deduce from \fer{ts.id}
\begin{align}
\label{est.time10}
\nor{\p_t u_{in}}_{L^{6/5}(0,T;(W^{1,6}(\O))')} \leq C.	
\end{align}

\no\textbf{Step 4. Identification of the limit. }
Since $u_{in}$ is a strong solution of \fer{eq.eq}-\fer{eq.id}, we have, for $\xi$ smooth,
\begin{align*}
  \int_0^T \p_t u_{in}(t,\cdot) \xi(t,\cdot)dt   & - \int_{Q_T} \Deltan p_i (\bu_n)  \xi = \int_{Q_T} f_i(\bu_n) \xi.
\end{align*}
Using the strong convergence \fer{lemma.conv1}, we easily justify the passing to the limit $n\to\infty$ for the reaction terms if $\xi\in L^{q/(q-2)}(Q_T)$. In view of \fer{est.time10}, the terms involving the time derivative are well defined and some subsequences (not relabeled) pass to their corresponding limits (weakly) if $\xi\in L^6(0,T;W^{1,6}(\O))$. Regarding  the diffusion term, we rewrite it as 
\begin{align*}
&\int_{Q_T} \Deltan p_i (\bu_n)  \xi  = \int_{Q_T} \int_\O J_n(x-y) \big(p_i (\bu_n(t,y))-p_i (\bu_n(t,x)) \big) dy~ \xi(t,x)dx dt \\
  &= \frac{1}{2} \int _{Q_T}  \int_\O J_n(x-y) \big(p_i(\bu_n(t,y))-p_i(\bu_n(t,x))\big) \big(\xi(t,y)-\xi(t,x)\big)dy dx dt \\
  &= \frac{C_1}{2} n^{2+N} \int _{Q_T}  \int_\O J(n(x-y)) \big(p_i(\bu_n(t,y))-p_i(\bu_n(t,x))\big) \big(\xi(t,y)-\xi(t,x)\big)dy dxdt \\
   &= \frac{C_1}{2} \int _{Q_T}\int _{\R^N} J(z) \chil_\O(x+\eps z)\frac{p_i(\bar \bu_n(t,x+\eps z))-p_i( \bu_n(t,x))}{\eps} \\
   & \hspace{3cm}\times \frac{\bar\xi(t,x+\eps z)-\xi(t,x)}{\eps}dz dxdt ,
\end{align*}
where $\eps=1/n$.
 We, again, examine the terms of $p_i(\bu)=c_i u_i +a_iu_i^2+u_iu_j$ separately. For the linear term, we have, using \fer{lemma.conv2},
\begin{align} 
\frac{C_1}{2} \int_{Q_T}\int _{\R^N}  & J(z) \chil_\O(x+\eps z)\frac{\bar u_{in}(t,x+\eps z)- u_{in}(t,x))}{\eps} \frac{\bar\xi(t,x+\eps z)-\xi(t,x)}{\eps}dz dx dt \nonumber  \\
&\to \frac{C_1}{2} \int_{Q_T}\int _{\R^N} J(z)~ z\cdot \grad u_i(t,x) ~ z\cdot \grad \xi(t,x) dz dx dt \nonumber\\
& = \int_0^T\int_\O \ba(\grad u_i(t,x)) \cdot \grad \xi(t,x) dxdt, \label{conv.f1}
\end{align}
where we defined, for $\bv\in\R^N$, 
\begin{align}
\label{def.a}
\ba_j(\bv)=\frac{C_1}{2}	\int _{\R^N}  J(z)~ z\cdot \bv ~ z_j dz = \bv, 
\end{align}
according to \cite[Lemma 6.16]{Andreu2010}.
Observe that the convergence \fer{conv.f1} may be extended, by density, to  $\xi\in L^{2}(0,T;H^1(\O))$. Similarly, for $u_iu_j$ we use, in addition to the weak convergence \fer{lemma.conv2}, the strong convergence deduced in Step 2. Splitting this term as in \fer{decomp}, we only have to examine, by symmetry, the following expression
\begin{align*}
I_{ij}=\frac{C_1}{2} \int_{Q_T}\int _{\R^N} &  J(z) \chil_\O(x+\eps z)\bar u_{jn}(t,x+\eps z)\frac{\bar u_{in}(t,x+\eps z)- u_{in}(t,x)}{\eps} \\  
&\times \frac{\bar\xi(t,x+\eps z)-\xi(t,x)}{\eps}dzdxdt.	
\end{align*}
Since $r'=2q/(q-2)$, we have
\begin{align*} 
I_{ij}=&\frac{C_1}{2} \int_{Q_T}\int _{\R^N}   \chil_\O(x+\eps z) J(z)^{1/q}\bar u_{jn}(t,x+\eps z)\\
& \times J(z)^{1/2}\frac{\bar u_{in}(t,x+\eps z)- u_{in}(t,x)}{\eps} \\  
&\times J(z)^{1/r'} \frac{\bar\xi(t,x+\eps z)-\xi(t,x)}{\eps}dzdxdt\\
&\to \frac{C_1}{2} \int_{Q_T} u_j(t,x) \int _{\R^N}   J(z) ~z\cdot \grad u_i(t,x) ~ z\cdot \grad \xi(t,x) dzdxdt  \\
&= \int_{Q_T} u_j \grad u_i \cdot \grad \xi ,
\end{align*}
where we used \fer{def.a} and \fer{lemma.conv2}. Observe that, in this case, we may extend the functional space of test functions, by density, to $\xi\in L^{r'}(0,T;W^{1,r'}(\O))$, and that, in relation to the weak convergence of the time derivatives, see \fer{est.time10}, we have $r'>6$.

Finally, a new and straightforward duality calculation shows that the initial data may be interpreted in the sense of \fer{idconv}. This finishes the proof of Theorem~\ref{th:convergence}.

\section{Proofs of the lemmas}\label{sec:lemmas}

\no\textbf{Proof of Lemma~\ref{lemma:dual}. }
The proof is similar to that of \cite[Lemma 3.8]{Andreu2010}. Fix $t_0>0$ and consider the Banach space $X_{t_0}=C([0,t_0];L^\infty(\O))$. Consider the operator 
\begin{align*}
\cT(w)(t,x)	= c(x)+\int_0^t  \Big(a(s,x)\int_\O \rho(x-y) (w(s,y)-w(s,x))dy + b(s,x)\Big)ds.
\end{align*}
To apply Banach's fixed point theorem we must check: (i) $\cT: X_{t_0}\to X_{t_0}$, and (ii) $\cT$ is contractive. We start with (i). For $0\leq t_1<t_2\leq t_0$, we have
\begin{align*}
\nor{\cT(w)(t_2,\cdot)-\cT(w)(t_1,\cdot)}_{L^\infty(\O)}	\leq &  k (t_2-t_1),
\end{align*}
where $k= \abs{\O}\nor{a}_{L^\infty(Q_{t_0})}\big(2 \nor{J}_{L^\infty(\R^N)} \nor{w}_{L^\infty(Q_{t_0})} +\nor{b}_{L^\infty(Q_{t_0})}\big)$. Similarly,
\begin{align*}
\nor{\cT(w)(t,\cdot)-c}_{L^\infty(\O)}	\leq	  k t.
\end{align*}
These two estimates give that $\cT(w) \in X_{T_0}$. To prove (ii), let $w,z\in X_{t_0}$. Then, for $t\in (0,t_0)$,
\begin{align*}
\nor{\cT(w)(t,\cdot)- \cT(z)(t,\cdot)}_{L^\infty(\O)}\leq 2 t \nor{a}_{L^\infty(Q_{t_0})} \nor{\rho}_{L^\infty(\R^N)}
\nor{w-z}_{C([0,T];L^\infty(\O))}.
\end{align*}
Thus, choosing $t_0 < (2 \nor{a}_{L^\infty(Q_T)} \nor{\rho}_{L^\infty(\R^N)})^{-1}$ we deduce that $\cT$ is a strict contraction. Banach's fixed point theorem allows to deduce the existence of a unique solution, $\phi_1$, of \fer{eq.test1}-\fer{eq.test2} for $t\in[0,t_0]$. Replacing \fer{eq.test2} by $\phi(t_0,x)=\phi_1(t_0,x)$ and the time interval $[0,t_0]$ by $[t_0,2t_0]$ we again deduce the existence of a unique solution, $\phi_2$, of \fer{eq.test1}-\fer{eq.test2} for $t\in[t_0,2t_0]$. Continuing this procedure we obtain a solution of \fer{eq.test1}-\fer{eq.test2} defined on $[0,T]$. 

Regarding the regularity of the solution, since $\phi \in C([0,T];L^\infty(\O))$ and $\rho \in L^\infty(\R^N)$, we deduce that $\rho*\phi(t,\cdot) \in L^\infty(\O)$ for all $t\in [0,T]$. Therefore $\Delta^{\!1,\rho} \phi=\rho*\phi -\phi\int_\O \rho(\cdot-y)dy \in  C([0,T];L^\infty(\O))$, and then, from the equation \fer{eq.test1}, we deduce that $\p_t \phi \in  L^\infty(Q_T)$. 

Finally, assume that $a,b,c$ are non-negative and suppose that $\phi$ is negative somewhere. Let $\xi(t,x)=\phi(t,x) +\epsilon t$, with $\epsilon>0$ small enough so that $\xi$ is negative somewhere. Let $(t_0,x_0)$ be a point where $\xi$ attains its negative minimum. Then, $t_0>0$ and
\begin{align*}
\p_t\xi(t_0,x_0)=& \p_t \phi(t_0,x_0)+\epsilon	 \\
&>  a(t_0,x_0)\int_\O \rho(x_0-y) (\phi(t_0,y)-\phi(t_0,x_0))dy + b(t_0,x_0)\\
&=  a(t_0,x_0)\int_\O \rho(x_0-y) (\xi(t_0,y)-\xi(t_0,x_0))dy + b(t_0,x_0)\geq 0,
\end{align*}
 which is a contradiction. Therefore, $\phi\geq 0$ in $Q_T$.
$\Box$

\medskip

\no\textbf{Proof of Lemma~\ref{lemma:compactness}. }
From Theorem~\ref{th.brezis}, \fer{BQ1} yields
\begin{align}
&\int_{Q_T}  \abs{f_n(t,x)-(f_n(t,\cdot)*\Phi_\delta)(x)}^2 dxdt\leq \frac{C}{n^2}.	\label{est.brezis2}
\end{align}
Estimate  \fer{BQ0} implies the existence of a subsequence $\{f_{n_k}\}\subset L^3(Q_T)$ and a function $f \in L^3(Q_T)$ such that $f_{n_k}\wto f$ weakly in $L^3(Q_T)$. Moreover, estimates \fer{BQ2} and \fer{est.brezis2} ensure, see \cite[Lemma 5.1]{Lions}, the convergence  $f_{n_k}^2 \to f^2$ in the sense of distributions in $Q_T$ (we take $g^n=h^n=f_{n}$ in \cite[Lemma 5.1]{Lions}). This is,
\begin{align*}
\int_{Q_T} f_{n_k}^2\zeta \to \int_{Q_T} f^2\zeta \qtext{for all }\zeta\in C_c^\infty(Q_T).	
\end{align*}
Since $f_{n_k} \wto f$ weakly in $L^3(Q_T)$, we have 
\begin{align*}
\int_{Q_T} f_{n_k}^2\zeta \leq \nor{f_{n_k}}_{L^3(Q_T)}^2\nor{\zeta}_{L^3(Q_T)}\leq C\nor{\zeta}_{L^3(Q_T)},
\end{align*}
so that, passing to a new subsequence if required, we have $f_{n_k}^2\wto f^2$ weakly in $L^{3/2}(Q_T)$, by density.  As a consequence, we deduce that 
\begin{align*}
\int_{Q_T} \abs{f_{n_k}-f}^2 = 	\int_{Q_T} f_{n_k}^2 -2 \int_{Q_T} f_{n_k}f + \int_{Q_T} f^2 \to 0,
\end{align*}
since $1 \in L^{3}(Q_T)$ and $f\in L^{3/2}(Q_T)$. Then, the uniform bound \fer{BQ0} allows to obtain \fer{lemma.conv1}. $\Box$

\medskip

\no\textbf{Proof of Lemma~\ref{w023}. } 
We prove the result for a general power $p$, with $1<p<\infty$. By density, it is enough to show that 
$ \nor{\Deltan \psi}_{L^p(\O)} \leq C \nor{\psi}_{W^{2,p}_0(\O)}$ for $\psi \in C_c^\infty(\O)$.

Considering the extension of $\psi$ to $\R^N$ and splitting $\R^N$ in terms of $\O$ and $\O^c$, we obtain, using the triangle inequality,
\begin{align}
\int_{\O}\Big\lvert \Deltan \psi(x)\Big\rvert^p dx \leq C \Big(\int_{\R^N} &\Big\lvert \int_{\R^N} J_n(x-y) (\bar \psi(y)-\bar\psi(x))dy\Big\rvert^pdx \nonumber
 \\
&+ \int_{\O} \Big\lvert \int_{\O^c} J_n(x-y) (\bar \psi(y)-\bar\psi(x))dy\Big\rvert^pdx \Big). \label{l3.aq}
\end{align}

\no\textbf{Step A.} We estimate the first term of the right hand side of \fer{l3.aq}. Using the changes $n(y-x) = z$, $n(y-x) = -z$, and setting $\eps=1/n$, we get 
\begin{align} 
\int_{\R^N} J_n(x-y)& (\bar \psi(y)-\bar\psi(x))dy\nonumber \\
&=  \frac{C_1}{2} \int _{\R^N}  J(z) \frac{\bar\psi(x+\eps z)-2\bar\psi(x)+\bar\psi(x-\eps z)}{\eps^2}dz  .\label{def.psi}
\end{align}
We define, for $s\in[0,1]$ and $\sigma\in[-1,1]$, the functions $v(s) = \bar\psi(x+s\eps z)+\bar\psi(x-s\eps z)$ and $w(\sigma) = \grad \bar\psi(x+\sigma s \eps z)$. We have 
\begin{align*}
	\bar\psi(x+\eps z)- & 2\bar\psi(x)+\bar\psi(x-\eps z) = \int_0 ^1 v'(s)ds \\
	&= \eps \int_0^1  z^T\cdot ( \grad \bar\psi(x+s\eps z) - \grad \bar\psi(x-s\eps z))ds\\
		&= \eps \int_0^1   \int_{-1}^1 z^T \cdot w'(\sigma) d \sigma ds = \eps^2  \int_0^1 s    \int_{-1}^1 z^T D^2\bar\psi(x+\sigma s \eps z) z ~d\sigma ds .
\end{align*}
Then, Jensen's inequality yields
\begin{align*}
	\Big\lvert \frac{\bar\psi(x+\eps z)-  2\bar\psi(x)+\bar\psi(x-\eps z)}{\eps^2}\Big\rvert^p \leq 2^{p-1} \int_0^1     \int_{-1}^1 \lvert z^T D^2\bar\psi(x+\sigma s \eps z) z \rvert^p ~d\sigma ds,
\end{align*}
and on noting that $\sigma s \eps z$ is independent of $x$, we deduce 
\begin{align*}
 \int_{\R^N} \lvert z^T D^2\bar\psi(x+\sigma s \eps z) z \rvert^p ~d\sigma ds \leq \abs{z}^{2p}\nor{\bar\psi}^p_{W^{2,p}(\R^N)}.
\end{align*}
Therefore,  taking into account that the integration in \fer{def.psi} may be limited to $z\in B_r$, see \fer{suppJ}, and  applying H\"olders inequality in \fer{def.psi}, we deduce
\begin{align}
 \int_{\R^N} \Big\lvert &\int_{\R^N} J_n(x-y) (\bar \psi(y)-\bar\psi(x))dy\Big\rvert^pdx \nonumber \\
 & \leq C 	\Big( \int_{B_r}J^{p'}\Big)^{p/p'} \int_{\R^N} \int_{B_r} \Big\lvert \frac{\bar\psi(x+\eps z)-2\bar\psi(x)+\bar\psi(x-\eps z)}{\eps^2} \Big\rvert^p dz dx \nonumber  \\ 
& \leq C \nor{\bar\psi}^p_{W^{2,p}(\R^N)}. \label{est.bieeeen}
\end{align}

\no\textbf{Step B.} We estimate the second term of the right hand side of \fer{l3.aq}.
The integration is performed in a band enclosing $\p\O$. We define the bounded sets
\begin{align*}
&D_\eps = \{x\in\O: x+\eps z\in\O^c, \text{ for }z\in B_r\} ,\\
&F_\eps=\{y\in\O^c: y=x+\eps z, \text{ for }x\in\O,~z\in B_r\}.
\end{align*}
Observe that if $x\in \O\backslash D_\eps$ and $y\in \O^c$ or if  $x\in \O$ and  $y\in \O^c\backslash F_\eps$ then $J_n(x-y)=0$, since $\eps=1/n$. 
Therefore,
\begin{align}
\int_{\O} \Big\lvert \int_{\O^c} J_n(x-y)& (\bar \psi(y)-\bar\psi(x))dy\Big\rvert^pdx=
\int_{\O} \abs{\psi(x)}^p\Big\lvert \int_{\O^c} J_n(x-y) dy\Big\rvert^pdx \nonumber \\
&= C_1^p n^{2p}\int_{D_\eps} \abs{\psi(x)}^p\Big\lvert \int_{F_\eps} \tilde J_n(x-y) dy\Big\rvert^pdx \nonumber \\
& \leq 
C_1^p n^{2p}\int_{D_\eps} \abs{\psi(x)}^p dx,\label{l3.eq2}
\end{align}
because $\nor{\tilde J}_{L^1(\R^N)}=1$, see \fer{def.jtilde}.

We recall here the uniform cone property, enjoyed by Lipschitz sets \cite[Definition~6.3]{Delfour}: For all $x\in\p\O$, there exist positive numbers $h,\o$ and $\rho$ such that for all $y\in B_\rho(x)\cap \overline{\O}$, we have  that the cone of vertex $y$, heigth $h$, and aperture $\o$, denoted by $C_y(h,\o)$, satisfies $C_y(h,\o)\subset \O$. Symmetrically,  for all $z\in B_\rho(x)\cap \O^c$, we have  $C_z(h,\o)\subset \interior(\O^c)$.

We claim that 
\begin{align*}
D_\eps\cup F_\eps = \{x\in\R^N: \dist(x,\p\O)<r\eps\}.	
\end{align*}
Let us prove, for instance, that $F_\eps=\{x\in\O^c: \dist(x,\p\O)<r\eps\}$. 
On one hand, suppose that $y\in F_\eps$ but $\dist(y,\p\O)\geq r\eps$.  Then $B_{r\eps}(y)\cap\O =\emptyset$ but $\abs{y-x}<r\eps $, which is a contradiction, since $x\in\O$. On the other hand, let $y\in\O^c$ and $x_0\in\p\O$ be such that $\abs{y-x_0}=\dist(y,\p\O)\leq\beta r\epsilon$, with $\beta<1$. Notice that $x_0$ does exist because $\p\O$ is closed. Then, due to the uniform cone property, there exists a cone $C_{x_0}(h,\o)\subset \O$, so that $\abs{y-x}\leq \abs{y-x_0}+\abs{x_0-x} < \beta r\eps+(1-\beta)r\eps$ for all $x\in C_{x_0}(h,\o) \cap B_{\rho_0}(x_0)$, for $\rho_0 = (1-\beta)r\eps$. Thus, $y=x+\eps z$ for some $z\in B_r$.
A similar proof stands for the identity $D_\eps=\{x\in\O: \dist(x,\p\O)<r\eps\}$.

 Consider the collection of open balls $\cB=\{B_{2r\eps}(x)\}_{x\in \p\O}$. It is clear that 
$\overline{D_\eps\cup F_\eps}$ is covered by $\cB$. Since $\overline{D_\eps\cup F_\eps}$ is closed and bounded, and hence compact, we may extract a finite collection  
$\cB_{K}=\{B_{2r\eps}(x_{k})\}_{k=1}^{K}$ covering $\overline{D_\eps\cup F_\eps}$. Moreover, applying Vitali's covering lemma (finite version, see \cite[Theorem~8.5]{Rudin}), we may extract a sub-collection of disjoint balls $\cB_{K'}=\{B_{2r\eps}(x_{k'})\}_{k'=1}^{K'}$ such that $\cB_{V}=\{B_{6r\eps}(x_{k'})\}_{k'=1}^{K'}$ covers $\overline{D_\eps\cup F_\eps}$. 

In the following, we remove the primes from the indices and introduce the notation $\alpha =6r$ and $B_{\alpha\eps}^{k} = B_{6r\eps}(x_{k})$. 
It is easy to check that the collection $\cB_V$ satisfies the following properties: (i) $\abs{B_{\alpha\eps}^k\cap F_\eps} >0$, for all $k=1,\ldots,K$, and (ii) For all $g\in L^1(\R^N)$, there exists a constant $C>0$ independent of $\eps$ such that 
\begin{align*}
\sum_{k=1}^K \int_{B_{\alpha\eps}^k} \abs{g(x)}dx \leq C \nor{g}_{L^1(\R^N)}.
\end{align*}
Property (i) is, again, a consequence of the uniform cone property, while property (ii) follows from the collection $\{B_{2r\eps}^k\}_{k=1}^K$ being disjoint and from  the finite number of balls of radius $2r\eps$ contained in a ball of radius $6r\eps$. In particular, notice that $\abs{\cup_{k=1}^K B^k_{6r\eps}}\leq 3^N \abs{\cup_{k=1}^K B^k_{2r\eps}}$, see \cite[Theorem~8.5]{Rudin}

Property (i) ensures that $\bar\psi$ vanishes in the positive measure set $B_{\alpha\eps}^k\cap F_\eps$, so that the Poincar\'e's inequality yields
\begin{align*}
\int_{D_\eps} \abs{\psi(x)}^p dx&= \int_{D_\eps\cup F_\eps} \abs{\bar\psi(x)}^p dx\leq \sum_{k=1}^K \int_{B_{\alpha\eps}^k} \abs{\bar\psi(x)}^p	dx \\
&\leq P_{\alpha\eps}^p \sum_{k=1}^K \int_{B_{\alpha\eps}^k} \abs{\grad \bar\psi(x)}^p	dx,
\end{align*}
where $P_{\alpha\eps}$ is the constant of Poincar\'e (for the $p$-Laplacian) corresponding to the open ball  $B_{\alpha\eps}$. According to \cite[Chapter~5, Theorem~2]{Evans}, $P_{\alpha\eps}= C \alpha\eps$, where $C$ only depends on $N$ and $p$. On noting that the function $g(x)=\abs{\grad \bar\psi(x)}$ vanishes in $B_{\alpha\eps}^k\cap F_\eps$ (because $\psi\in C^\infty_c(\O)$ or, alternatively, $\psi\in W^{2,p}_0(\O)$), we may use again the Poincar\'e's inequality to obtain, on noting property (ii), 
\begin{align*}
\int_{D_\eps} \abs{\psi(x)}^p dx&\leq C P_{\alpha\eps}^{2p} \sum_{k=1}^K \int_{B_{\alpha\eps}^k} \abs{D^2 \bar\psi(x)}^p	dx \leq C\eps^{2p} \nor{\bar\psi}_{W^{2,p}(\R^N)}^p,
\end{align*}
 where we used that
 \begin{align*}
 \abs{\grad g} = \Big( \sum_{j=1}^N 	\Big\lvert \frac{\grad\bar\psi}{\abs{\grad\bar\psi}}\cdot D_j^2\bar\psi \Big\rvert^2\Big)^{1/2}.
 \end{align*}
Returning to \fer{l3.eq2} and noting that $\eps=1/n$,   we obtain 
\begin{align}
\int_{\O} \Big\lvert \int_{\O^c} J_n(x-y)& (\bar \psi(y)-\bar\psi(x))dy\Big\rvert^pdx \leq C\nor{\bar\psi}_{W^{2,p}(\R^N)}^p \label{mejoooor}
\end{align}

\no\textbf{Step C.} We finish the proof by replacing \fer{est.bieeeen} and \fer{mejoooor} in \fer{l3.aq} and recalling that $\nor{\bar\psi}_{W^{2,p}(\R^N)} = \nor{\psi}_{W^{2,p}_0(\O)}$, see \cite[Lemma~3.22]{Adams}. $\Box$

\begin{remark}
Taking the limit $n\to\infty$ in \fer{def.psi}, we obtain, for $\psi:\R^N\to\R$ smooth
\begin{align*} 
\Deltan \psi(x) \to & \frac{C_1}{2} \int _{\R^N}  J(z) z^T D^2\psi(x) z dz  .
\end{align*}
Since the Hessian of $\psi$,  $D^2\psi$, is symmetric, there exists an orthogonal (rotation) matrix, $R(x)$, and a diagonal matrix $Q(x)$ such that $D^2\psi(x) = R(x)^T Q(x) R(x)$, with $\det(R(x))=1$. Thus
\begin{align*}
\int _{\R^N}  J(z) z^T D^2\psi(x) z dz & =   	\int _{\R^N}  J(z) (R(x)z)^T Q(x) R(x)z dz \\
& =   \int _{\R^N}  J(Q^{-1}(x)y) y^T Q(x) y dy = \int _{\R^N}  J(y) y^T Q(x) y dy ,
\end{align*}
since $J$ is radial. Thus, 
\begin{align*} 
\frac{C_1}{2}\int _{\R^N}  J(z) z^T D^2\psi(x) z dz  = \sum_{i=1}^N Q_{ii}(x) \frac{C_1}{2}\int _{\R^N}  J(y) y_i^2 dy = \tr(Q(x)),
\end{align*}
where we used the normalization condition \fer{def.jeps}. Finally, since the trace is invariant under diagonalization, we deduce $\Deltan \psi(x) \to \Delta \psi(x)$ uniformly in $\R^N$.
	
\end{remark}

\no\textbf{Acknowledgements. } 
We acknowledge the financial support from the Spanish Ministerio de Ciencia e Innovaci\'on  project PID2020-116287GB-I00.


\begin{thebibliography}{9}

\bibitem{Adams}
Adams, R.A.:
Sobolev spaces.
Academic Press (1975)

\bibitem{Amann1989}
H.~Amann.:
  Dynamic theory of quasilinear parabolic systems. {III}. {G}lobal
  existence.
  Math. Z. 202, 219--250 (1989)

\bibitem{Andreu2010}
Andreu-Vaillo, F.,  Maz{\'o}n, J.M.,  Rossi, J.D.,  Toledo-Melero, J.J.:
  Nonlocal diffusion problems.
  American Mathematical Society (2010)

\bibitem{Bourgain2001} 
Bourgain, J. Brezis, H. Mironescu, P.:
  Another look at Sobolev spaces. 
  In: Menaldi, J.L., (eds.) Optimal Control and Partial Differential Equations,  439--455, IOS Press (2001)

\bibitem{Chen2004}
Chen, L., J{\"u}ngel, A.:
  Analysis of a multidimensional parabolic population model with strong
  cross-diffusion.
  SIAM J. Math. Anal. 36(1), 301--322 (2004)


\bibitem{Delfour}
Delfour, M.C., Zol\'esio, J.P.: Shapes and geometries. SIAM (2011)

\bibitem{Delia2020A}
D'Elia, M., Du, Q., Glusa, C., Gunzburger, M., Tian, X., Zhou, Z.: 
Numerical methods for nonlocal and fractional models. 
Acta Numer. 29, 1--124 (2020)

\bibitem{Delia2020B}
D'Elia, M., Gunzburger, M., Vollmann, C.: 
A cookbook for finite element methods for nonlocal problems, including quadrature rules and approximate Euclidean balls. 
Math. Models Methods Appl. Sci. 31(08), 1505--1567 (2021)







\bibitem{Desvillettes2014}
Desvillettes, L., Lepoutre, T., Moussa, A.:
  Entropy, duality, and cross diffusion.
  SIAM J. Math. Anal. 46(1), 820--853 (2014)

\bibitem{Du2019}
Du, Q., Yin, X.:
A conforming DG method for linear nonlocal models with integrable kernels.
J. Sci. Comput. 80,  1913--1935 (2019)


\bibitem{Evans}
Evans, L.C.:
Partial differential equations.
American Mathematical Society (1998)

\bibitem{Galiano2019a}
Galiano, G.:
  Well-posedness of an evolution problem with nonlocal diffusion.
  Nonlinear Anal. Real World Appl. 45, 170--185 (2019)

\bibitem{Galiano2021}
Galiano, G.:
Error analysis of some nonlocal diffusion discretization schemes.
Comput. Math. Appl. To appear.

\bibitem{Galiano2003}
Galiano, G., Garz{\'o}n, M.L., J{\"u}ngel, A.:
  Semi-discretization in time and numerical convergence of solutions of
  a nonlinear cross-diffusion population model.
  Numer. Math. 93(4), 655--673 (2003)

\bibitem{Galiano2014}
Galiano, G., Selgas, V.:
  On a cross-diffusion segregation problem arising from a model of interacting particles.
  Nonlinear Anal. Real World Appl. 18, 34--49 (2014)


\bibitem{Galiano2019b}
Galiano, G., Velasco, J.:
  Well-posedness of a cross-diffusion population model with nonlocal diffusion.
  SIAM J. Math. Anal. 51(4), 2884--2902 (2019)



\bibitem{Jungel2015}
J{\"u}ngel, A.:
The boundedness-by-entropy method for cross-diffusion systems. 
Nonlinearity 28, 1963–-2001 (2015)


\bibitem{Lions} 
Lions, P.L.:
Mathematical Topics in Fluid Mechanics: Volume 2: Compressible Models.
  Oxford Lecture Series in Mathematics and Its Applications, Clarendon Press (1996)


\bibitem{Perez2011}
P\'erez-Llanos, M., Rossi, J.D.:
Numerical approximations for a nonlocal evolution equation,
SIAM J. Numer. Anal. 49(5),  2103--2123 (2011)


\bibitem{Pierre1997}
Pierre, M., Schmitt, D.:
  Blowup in reaction-diffusion systems with dissipation of mass.
  SIAM J. Math. Anal. 28, 259--269 (1997)

\bibitem{Rudin}
Rudin, W.: Real and complex analysis.
McGraw-Hill (1970)


\bibitem{Shigesada1979}
Shigesada, N., Kawasaki, K. Teramoto, E.:
  Spatial segregation of interacting species.
  J. Theoret. Biol. 79(1), 83--99 (1979)


\end{thebibliography}
\end{document}